\documentclass[preprint,review,12pt,sort&compress]{elsarticle}
\pdfoutput=1
\usepackage[margin=1in]{geometry}
\usepackage{multirow}
\usepackage{amsmath}
\usepackage{cases}
\usepackage{amssymb}
\usepackage{gensymb}
\usepackage{bm}
\usepackage{enumitem}
\usepackage{lipsum}
\usepackage{makecell}
\usepackage{soul}

\usepackage{amssymb}

\usepackage{float}

\usepackage{xcolor}

\journal{Materials and Design}

\begin{document}

\begin{frontmatter}


\cortext[cor1]{Corresponding Author}

\title{Graph-theoretic estimation of reconfigurability in origami-based metamaterials}


\author[address1]{Koshiro Yamaguchi}
\author[address1,address2]{Hiromi Yasuda}
\author[address3]{Kosei Tsujikawa}
\author[address4]{Takahiro Kunimine}
\author[address1]{Jinkyu Yang\corref{cor1}}
\ead{jkyang@aa.washington.edu}

\address[address1]{William E. Boeing Department of Aeronautics and Astronautics, University of Washington, Seattle, Washington 98195-2400, USA}
\address[address2]{Department of Mechanical Engineering and Applied Mechanics, University of Pennsylvania, Philadelphia, PA 19104, USA}
\address[address3]{Division of Mechanical Science and Engineering, Graduate School of Natural Science and Technology, Kanazawa University, Kanazawa, Ishikawa, 920-1192, Japan}
\address[address4]{Faculty of Mechanical Engineering, Institute of Science and Engineering, Kanazawa University, Kanazawa, Ishikawa, 920-1192, Japan}

\begin{abstract}
\linespread{1}\selectfont{
Origami-based mechanical metamaterials have recently received significant scientific interest due to their versatile and reconfigurable architectures. However, it is often challenging to account for all possible geometrical configurations of the origami assembly when each origami cell can take multiple phases. Here, we investigate the reconfigurability of a tessellation of origami-based cellular structures composed of bellows-like unit cells, specifically Tachi-Miura Polyhedron (TMP). One of the unique features of the TMP is that a single cell can take four different phases in a rigid foldable manner. Therefore, the TMP tessellation can achieve various shapes out of one originally given assembly. To assess the geometrical validity of the astronomical number of origami phase combinations, we build a graph-theoretical framework to describe the connectivity of unit cells and to analyze the reconfigurability of the tessellations. Our approach can pave the way to develop a systematic computational tool to design origami-based mechanical metamaterials with tailored properties.
}
\end{abstract}


\end{frontmatter}


\section{Introduction}
\label{intro}
Mechanical metamaterials can have a wide range of mechanical properties by leveraging the design freedom in their architectures~\cite{Babaee2015, Bilal2017, Wu2018}. Origami has helped this emerging platform to add the tunability of its shapes and behaviors after fabrications \cite{Schenk2013, Felton2014}. Recently, an extensive amount of research on origami-based metamaterials has been conducted to show their various characteristics in both mathematical and mechanical aspects \cite{tachi_introduction_2019, Rus2018}. Particular interest has been placed on the realization of reconfigurable mechanical metamaterials via the concept of origami, kirigami, or equivalent  architectures~\cite{Yang2020, overvelde_three-dimensional_2016,overvelde_rational_2017, Filipov2015,yasuda_origamibased_2019,Fang2018}. While most research focuses on the homogeneous arrangements of origami units and their intrinsic properties, the cross-linkage of heterogeneous origami cells and their effective global characteristics have been relatively unexplored. 

To address this connectivity and reconfigurability problem of heterogeneous origami architectures, Tachi-Miura Polyhedron (TMP) can serve as an ideal platform~\cite{Tachi2012, yasuda_origamibased_2019}. TMP is one of the origami-based mechanical metamaterials with bellow-like 3D unit cells constructed from 2D Miura-folding sheets. This TMP-based metamaterial features several mechanical characteristics, such as flat and rigid foldability, negative Poisson's ratio, multistability, and load-bearing capabilities~\cite{yasuda_reentrant_2015, yasuda_origamibased_2019}. Notably, each TMP cell can take four different states (i.e., phases) while maintaining rigid foldability. Thus, when we combine them into a tessellation, we can construct a 2D array of heterogeneous origami cells. However, we cannot arbitrarily assign one of four states to each unit cell randomly, because it can result in a geometrically invalid configuration. Therefore, we have to find a method to discover geometrically valid tessellations efficiently. If we follow the brute-force method, the number of possible configurations that need to be assessed becomes $4^n$, where $n$ is the total number of the TMP cells. This exponential increase of the search space is often called the ``curse of dimensionality", also known as the combinatorial explosion. This leads the search for the geometrically valid configuration impossible for a large tessellation. 

To overcome this hurdle, we can consult the graph-theoretical method to describe and analyze the TMP-based tessellations. Graph theory is a field in discrete mathematics, which is used to represent pairwise relations between objects~\cite{Mesbahi2010}. In this mathematical schematization, the objects are referred to as vertices, and the connected pairs of the objects are called edges. In diagrammatic ways, graphs are often represented with dots (or circles) for vertices and lines for edges. Generally, graphs are classified into two types: undirected graphs and directed graphs. The former have edges without direction (two-way relationship), whereas the latter has edges with direction (one-way relationship). This simple-yet-efficient mathematical structure can be utilized to model various types of connections in general. Therefore, it has been applied to numerous systems, such as physical~\cite{Deo1974}, biological~\cite{Canutescu2003}, and information~\cite{Xu2018} ones. However, this versatile feature of graph theory to model networks has not been fully leveraged to model and analyze tessellations of mechanical metamaterials composed of unit cells, to the best of the authors' knowledge. 

In this study, we employ the idea of enumerating the connection of the system composed of the TMP unit cells, and we search the unique configurations and count the number of them to assess the reconfigurability of the origami system. The results show the exact number of unique configurations even for a large tessellation, which could not have been calculated easily by the conventional brute-force approach. We also report the growing pattern of the number of valid configurations as the number of TMP cells increases in horizontal and vertical directions. Among those results, we find heterogeneous configurations composed of different states of unit cells. Unlike normal configurations that are homogenous and have less variety of unit cells, such heterogeneous configurations are hard to be discovered by human heuristics. These results suggest that the morphology of the reconfigurable mechanical metamaterials can meet the engineering requirements that are evolving and adapting to external environments.

\section{Methods and materials}
\subsection{TMP tessellation and its transformation}
\begin{figure}[H]
\centering
\includegraphics{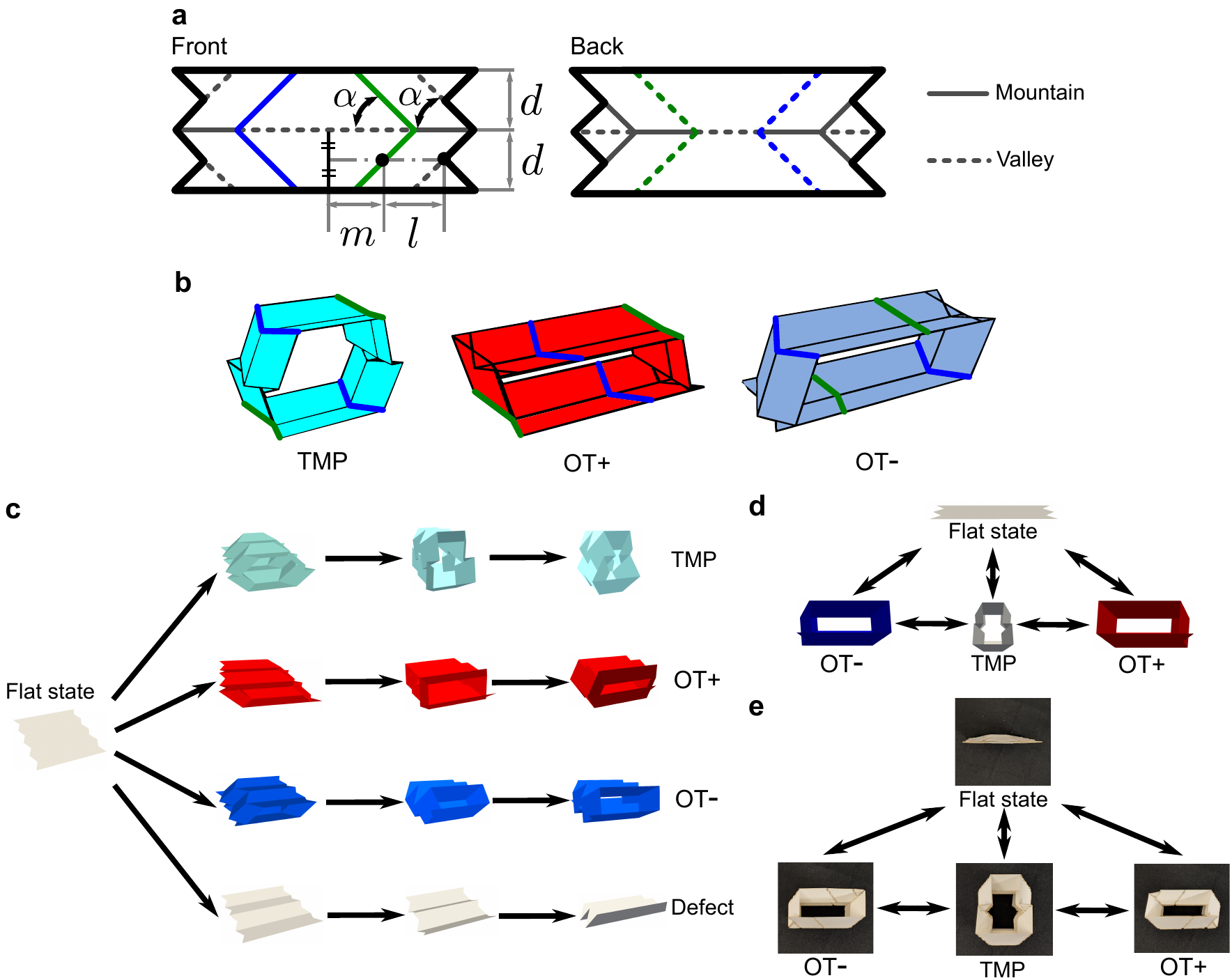}
\caption{Folding behavior of a single Tachi-Miura Polyhedron (TMP) unit cell. \textbf{a}, crease patterns and geometrical parameters of two flat sheets composing the TMP. \textbf{b}, the folding angles of blue and green creases decide the reconfigurable states of a unit cell. When blue (green) creases are folded flat, the unit cell takes the OT$+$ (OT$-$) phase. \textbf{c}, 3D rendered images of the transition from the flat state four reconfigurable states. \textbf{d}, 3D rendered images of a TMP unit cell with the initial flat state and the three reconfigurable states of TMP, OT$+$, and OT$-$. \textbf{e}, images of paper prototypes with corresponding configurations to \textbf{d}.}
\label{fig:morphing_scheme}
\end{figure}

It is well known that TMP can take finite volume or be collapsed to zero volume during its folding process \cite{yasuda_reentrant_2015} (see ``TMP" and ``flat" state in Fig.~\ref{fig:morphing_scheme}). In this study, we also consider the other two geometries, so-called origami tube minus (OT$-$, which is a left-parallelepiped tube) and origami tube plus (OT$+$, a right-parallelepiped tube), which can be generated from the flat state without deviating from the rigid-foldable assumption (see Supplementary Movie 1 for the demonstration of a paper prototype, and Supplementary Movie 2 for the 3D rendered transformation). Thus, in this study, we consider four possible states of unit cells: TMP, OT$-$, OT$+$, and Defect. Here, defect mode means the flat state, which does not occupy any volume. 

Fig.~\ref{fig:morphing_scheme}(a) explains the components of TMPs with geometrical parameters $(l,m,d,\alpha)$. In Fig.~\ref{fig:morphing_scheme}(b), we see that the folding angle of blue and green creases decides the reconfigurable states of a unit cell.  If either blue or green creases are folded flat, unit cells become OT$+$ or OT$-$, respectively. Fig.~\ref{fig:morphing_scheme}(c) shows the transition of the unit cell from the flat state to four reconfigurable states of TMP, OT$+$, OT$-$, and Defect. Figs.~\ref{fig:morphing_scheme}(d) and (e) depict those states and transformations between them for the 3D-rendered model and paper prototype, respectively. 

\begin{figure}[H]
\centering
\includegraphics{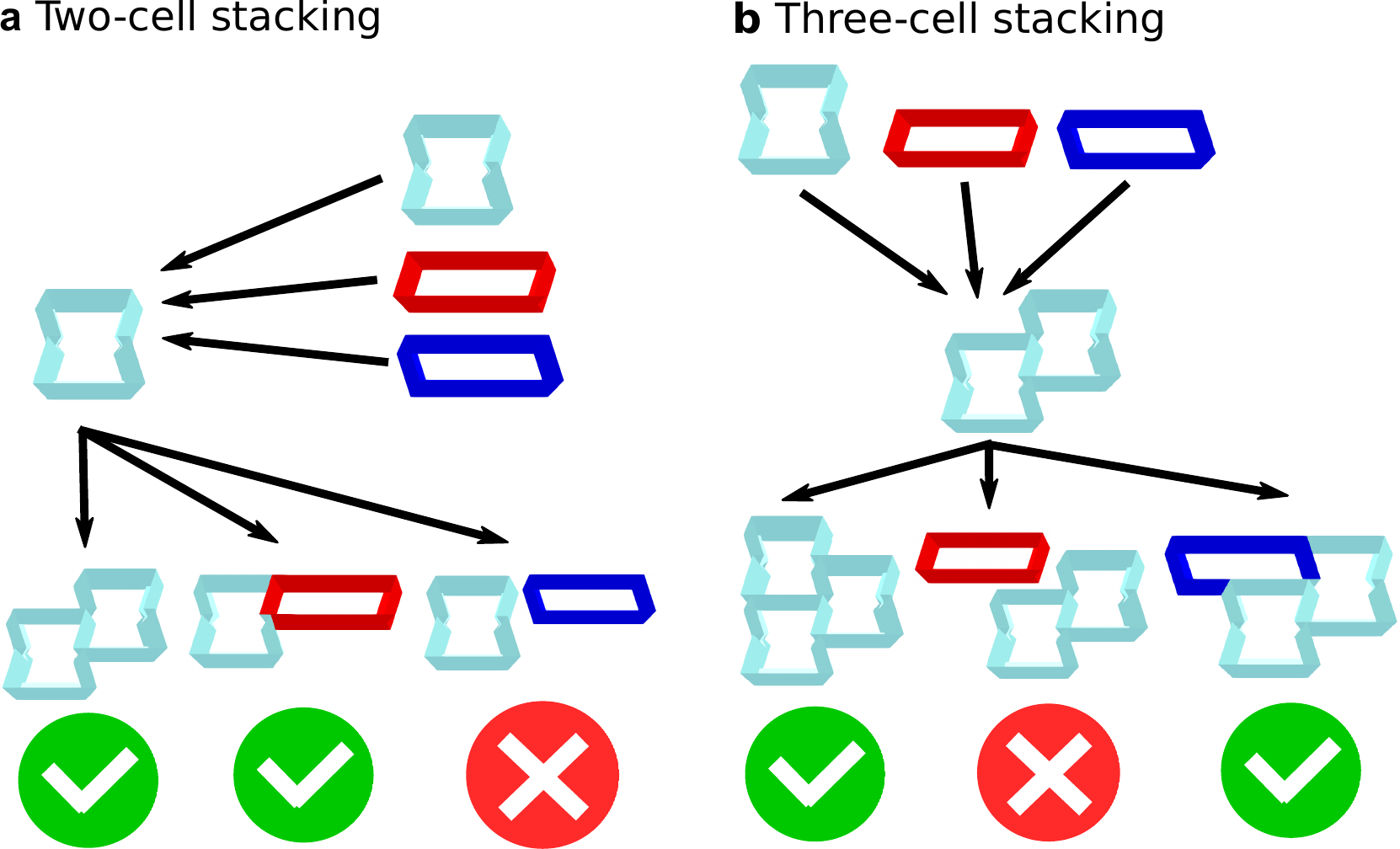}
\caption{Examples of how each tessellation becomes valid or invalid. Cases of \textbf{a}, two-cell tessellations and \textbf{b}, three-cell tessellations. In \textbf{a}, one unit cell (TMP, OT$+$, or OT$-$) is attached to one TMP cell. The attachment of TMP and OT$+$ results in the valid configuration. However, we can see that OT$-$ tube does not fit with the upper-right side of the TMP cell. Likewise, cases of three-cell stacking in \textbf{b} shows that the attachment of TMP and OT$-$ to a two-TMP stacking brings the valid configuration, whereas the attachment of OT$+$ tube makes the tessellation invalid. }
\label{fig:validorinvalid}
\end{figure}

The TMP cells can further be assembled into a tessellation in both horizontal and vertical manner as shown in Fig.~\ref{fig:tessellation}(a), and it can have multiple cells with different states (see Supplementary Movie 3 for a transformation of a tessellation). For example, Figs.~\ref{fig:tessellation}(b) and (c) show some configurations of 3D-rendered and paper 3-by-3 stacking with TMP, OT$-$, and OT$+$, which are represented by gray, blue, and red colors, respectively. Including the defect mode, we have four different states that can be assigned for each unit. Therefore, we can consider $4^n$ patterns of configurations within a $n$-cell tessellation. However, there exists a kinematic relationship between the adjoined TMP cells, and a large portion of the $4^n$ patterns are not valid. Fig.~\ref{fig:validorinvalid} illustrates this problem by using two-cell and three-cell tessellations. By simulating the stacking of unit cells, we can tell that some configurations are invalid due to the unmatching of the geometrical boundaries of the unit cells. 

While we can achieve various configurations and shapes with this reconfigurable origami tessellation, it is hard to know the actual number of total configurations within a given size of tessellation. Considering all cases (e.g., $4^9$ = 262,144 for 3-by-3 TMP tessellations) and sorting out what are valid configurations one by one is computationally expensive and even impossible for a large tessellation using conventional computational power. To overcome this hurdle, we introduce the idea of a graph-theoretical framework to discover the possible configuration efficiently. With the aid of the graph representation of a reconfigurable origami tessellation and adjacency matrix obtained from the graph information, we are able to calculate the total number of valid configurations systematically. In the following sections, we introduce the graph representation and adjacency matrices to assess the combinatorial problem of the TMPs.

\subsection{Graph representation of TMP tessellation}
Recently, several researchers have tackled tessellation-related problems for materials sciences with the aid of graph-theoretic techniques. Arkus et al. utilized graph-theoretic enumerations of adjacency matrices and distance matrices to find the finite sphere packings \cite{Arkus2011}. Vlassis et al. applied graphs and graph convolutional deep neural network for anisotropic hyperelasticity of polycrystal materials \cite{Vlassis2020}. Herein, we further leverage this technology to investigate the reconfigurability of origami tessellations. 

\begin{figure}[H]
\centering
\includegraphics{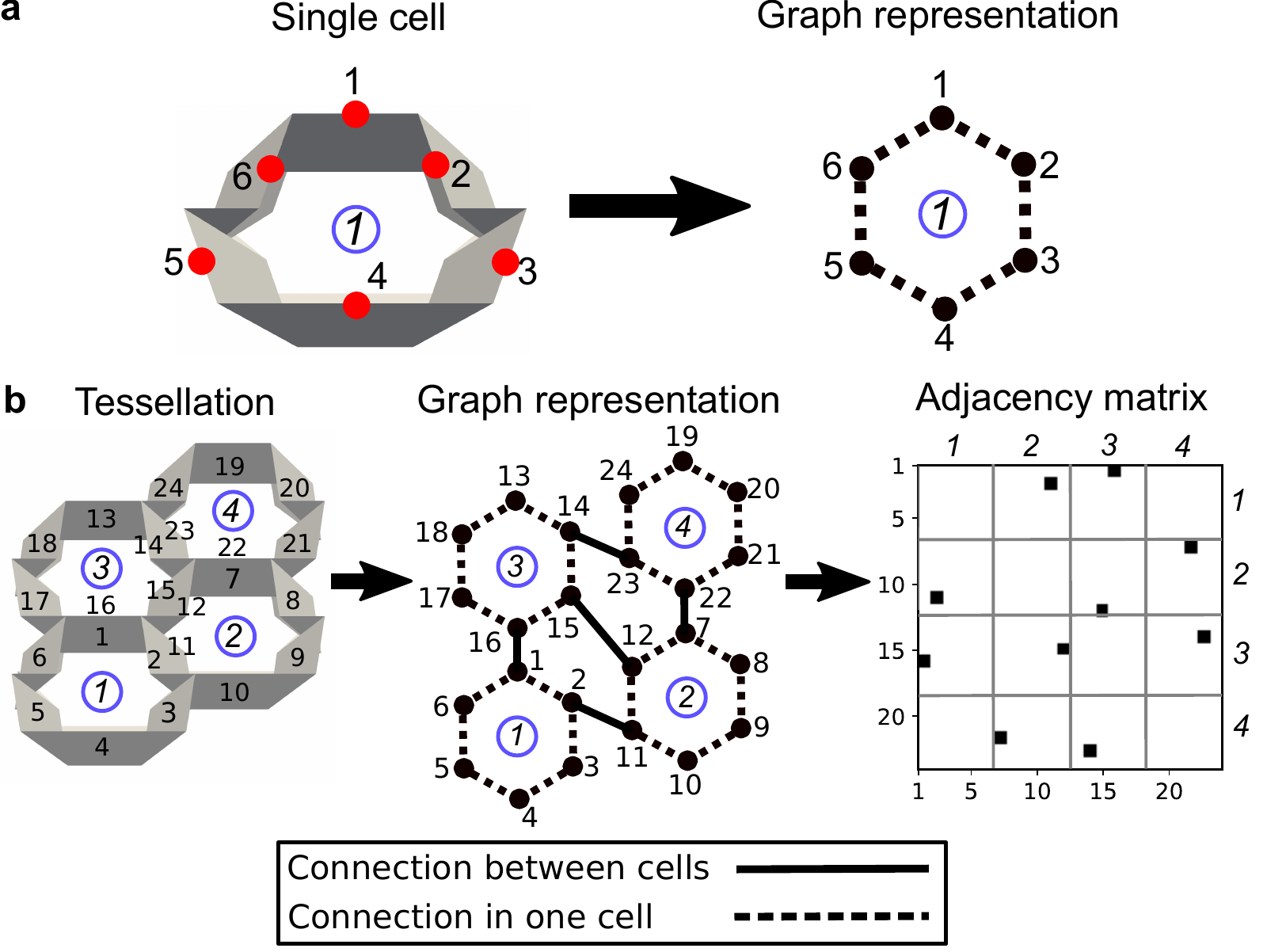}
\caption{Examples of the graph representation for \textbf{a}, a single unit cell and \textbf{b}, a 2-by-2 tessellation. Red dots in \textbf{a} represent the midpoints of the major side edges of TMPs. In a graph representation, those red dots are considered as nodes. Dashed lines in \textbf{a} and \textbf{b} mean the connection of the nodes within a unit cell. Solid lines in \textbf{b} means the connection of the nodes between different unit cells. Italic numbers in blue circles denote the TMP unit cells, and upright numbers denote the middle points in physical TMPs and corresponding nodes in graph. In the illustration of the adjacency matrix, left and bottom labels shows the node numbers as represented in the graph. Likewise, right and top labels in italic numbers denote the TMP unit cells. Gray lines in adjacency matrix show which part of the matrix corresponds to the TMP unit cell numbers. 
Dots in this adjacency matrix indicate the connected nodes, i.e., adjoined faces in the TMP tessellation.} 
\label{fig:graphrepresentation}
\end{figure}

As a first step, we attempt to translate physical TMP tessellations to graph representations. Fig.~\ref{fig:graphrepresentation} shows how we define an undirected graph $G=(V,E)$ from a single unit cell in Fig.~\ref{fig:graphrepresentation}(a) and a tessellation in Fig.~\ref{fig:graphrepresentation}(b), respectively. Here, vertices $V$ represents the six major faces of a unit cell, and the edge set $E$ consists of unordered pairs of vertices. Here, the number of vertices is $|V|=6n$, where $n$ is the number of unit cells in a tessellation. The edge set counts the connections of faces both within one unit cell and between other cells (see dashed and solid lines in Fig.~\ref{fig:graphrepresentation}). Then, we also define a subgraph $G'=(V',E')$ such that $V'\subseteq V$ and $E' \subseteq E$, and $V'$ and $E'$ represents the interconnections among the unit cells. For example, in Fig.~\ref{fig:graphrepresentation}(b) where we consider a 2 by 2 tessellation, we can see that $V'=\{1,2,7,11,12,14,15,16,22,23\}$ and $E'=\{(1,16),(2,11),(7,22),(12,15),(14,23)\}$.

Based on this graph representation, we utilize graph-theoretic enumeration of adjacency matrices. We introduce those definitions in the following.
A tessellation of $n$ unit cells can be described by an $|V|\times |V|$ adjacency matrix $\bm{A}$. This matrix is defined with the graph $G$ and the subgraph $G'$.
Therefore it details which vertices (faces in physical unit cells) are in contact: $\bm{A}_{ij}=1$ if and only if $(u_i,u_j) \in E'$ and otherwise $\bm{A}_{ij}=0$. This means that $\bm{A}_{ij}=1$ if $i$th and $j$th faces are attached to each other, and $\bm{A}_{ij}=0$ if they are not. It should be noted that each $x$-by-$y$ tessellation has only one unique adjacency matrix -- regardless of the phases of unit cells --  because once we define the size of the tessellation and its arrangement, we also know which physical faces of unit cells to be connected to others. The illustration of the adjacency matrix is also shown in Fig.~\ref{fig:graphrepresentation}(b). Here, the dots denote the edge set $E'$, which represents the adjoining of TMP faces.

\subsection{Combinatorial search for valid configurations} 
To analyze the reconfigurability of the TMP tessellations, we utilize the aforementioned graph representation and algorithm. The reconfigurability can be quantified by the number of unique configurations that the tessellation can achieve. The search of valid configurations consists of three steps: (i) definition of the size of tessellation, (ii) graph representation and graph-theoretic enumeration of the tessellation, and (iii) combinatorial search for valid configurations. First, we define the size and shape of the tessellation. Our graph-theoretic method is versatile and applicable to various types of tessellations. However, in order to make comparisons and measure the effect of the unit cell number on the reconfigurability, we use a rectangular arrangement of $X$-by-$Y$ tessellation, where $X$ and $Y$ are the number of unit cells in horizontal and vertical directions, respectively.

Secondly, we represent the tessellation via graph and enumerate the information about the connection of the unit cells using an adjacency matrix. Again, the adjacency matrix is unique for each tessellation and immutable in the process of tessellation search. After establishing the information about tessellations, we move to the final step, which is the combinatorial search for valid configurations. In this step, we utilize the complete list of 3-cell configurations that are known a priori. It has 32 configurations as shown in Supplementary Figure S2, and we combine those 3-cell configurations to find configurations of larger tessellations.

\begin{figure}[H]
\centering
\includegraphics{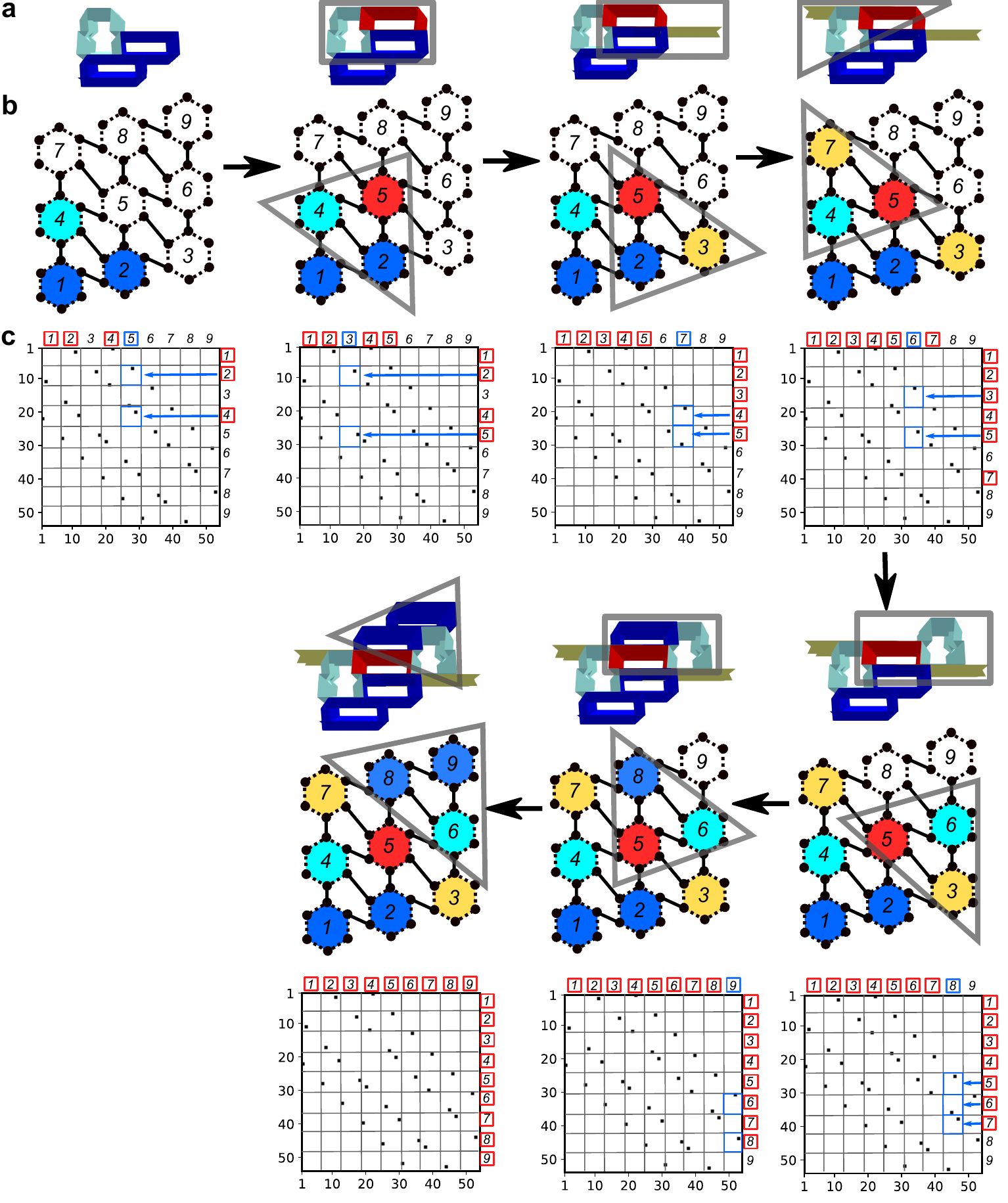}
\caption{An example of discovering a valid configuration for the 3-by-3 tessellation with logistic illustrations of the computation processes using adjacency matrices. In \textbf{a} and \textbf{b}, cyan, red, blue, and yellow colors in the graph representation correspond to the TMP, OT$+$, OT$-$, and defect states of the unit cell, respectively. In \textbf{a}, 3d-rendered images represent the graphical process of building tessellations. Likewise, in \textbf{b}, the pictures of graphs with colored parts show the schematic process of the tessellation search. Three unit cells enclosed by gray lines in \textbf{a} and \textbf{b} indicates where the 3-cell configurations are referenced and adopted. In \textbf{c}, adjacency matrices are illustrated in the same way as Fig.\ref{fig:graphrepresentation}(b). Italic numbers with red boxes show that unit cells have assignments of configurations, whereas italic numbers with no boxes show that those unit cells do not have assignments yet. Blue arrows and boxes in the matrix show the process of searching an adjacent unit cell to assign the proper unit cell configuration.} 
\label{fig:searchprocess}
\end{figure}

The schematic illustration of the process of this combinatorial search is shown in Fig.~\ref{fig:searchprocess}. The top row shows the graphical evolution of the TMP tessellation, the middle row shows their graph representation, and the bottom panel explains how we use adjacency matrices for this search process. Italic numbers with red boxes show that unit cells have assignments of configurations, whereas italic numbers with no boxes show that those unit cells do not have assignments yet. By using this assignment information and adjacency matrix, we detect an unassigned unit cell that is adjacent to the assigned ones. Blue arrows and boxes in the matrix show the process of searching for an adjacent unit cell to assign. To begin the search, we allocate configurations to three cells at the lower-left end. Then, based on the graph structure, we find an adjacent cell and assign a configuration to the adjacent one by looking up the list of 3-cell configurations (see gray triangles in Fig.~\ref{fig:searchprocess}b). We repeat this process until all unit cells have configurations. Also, in this process, we check every one of the configurations in the list. Therefore, at the end of the process, we can count all valid configurations out of a given size of tessellation.

\section{Results and discussion}
\subsection{Estimation and growth of reconfigurability}
\begin{figure}[H]
\centering
\includegraphics{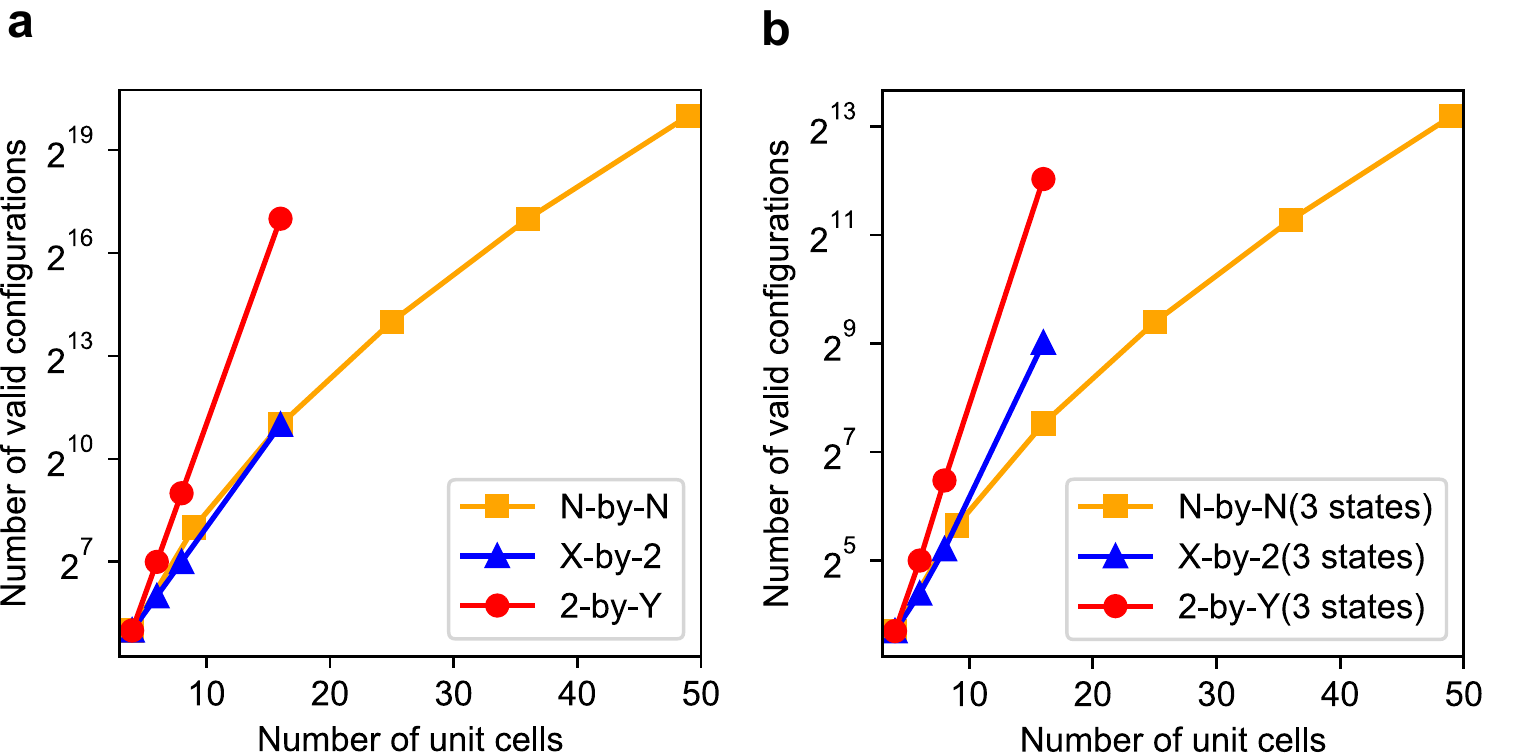}
\caption{Number of valid configuration for each tessellation size. \textbf{a}, 14 results of counting valid configurations with four configurations (TMP, OT$+$, OT$-$, and Defect) based on the categorization of N-by-N, X-by-2, and 2-by-Y tessellations where $X$, $Y$, and $N$ are the number of horizontal and vertical stackings of the tessellation. \textbf{b}, 14 results of counting valid configurations with three configurations (OT$+$/OT$-$/Flat, TMP/OT$+$/Flat, TMP/OT$-$/Flat, and TMP/OT$+$/OT$-$). Each three-state configuration has the same number of valid configurations. Other results and exact numbers are shown in Table \ref{tab:Combinations}.} 
\label{fig:numberofconfig}
\end{figure}

We start discussing the number and growth of reconfigurability by changing the size of tessellations. The result of the searching for the number of valid configurations with 14 different tessellation sizes is summarized in Fig.~\ref{fig:numberofconfig} for selected results and Table.~\ref{tab:Combinations} for all results. To verify the result from our graph-theoretical method, we also implement a brute-force method for smaller sizes of tessellations in Supplementary Note A. We start to examine the unique configurations with 2-by-2 tessellations with the total number of four unit cells. We then enlarge the size of the tessellation by adding a row or column of unit cells to 49 unit cells in total. If we follow the brute-force method for this number of unit cells, we have to investigate $4^{49} \approx 10^{29}$ candidates of configurations that are almost impossible to examine even with the aid of high-performance computers. Therefore, by utilizing our method, we achieve a substantial amount of efficiency in estimating the reconfigurability of the origami-based metamaterials. 

\begin{figure}[H]
\centering
\includegraphics{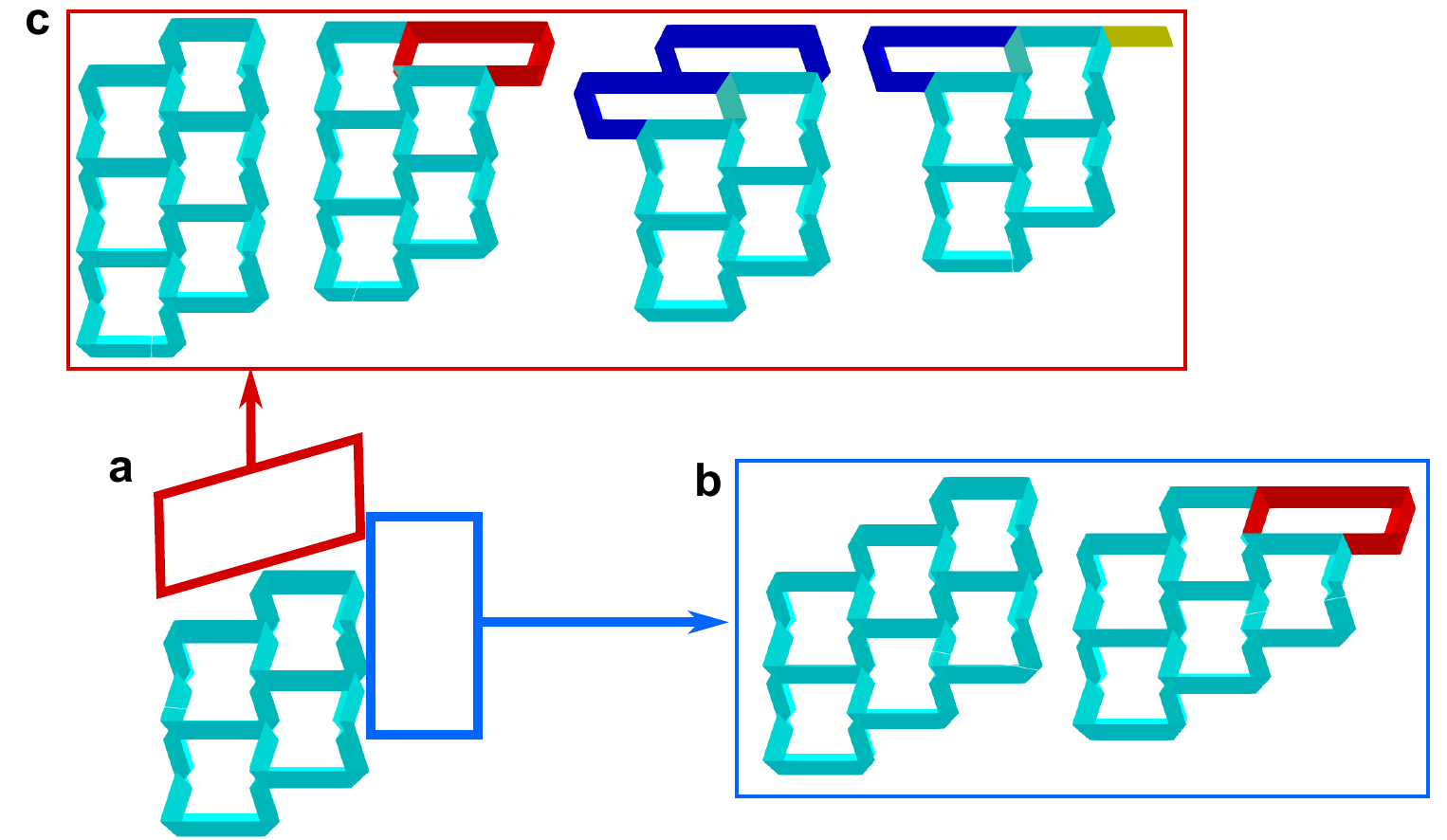}
\caption{\textbf{a}, an image of a 2-by-2 tessellation with four TMPs. Red box shows the position where the unit cells are attached to build 2-by-3 tessellations. Likewise, blue box shows the position for 3-by-2 tessellations. \textbf{b}, 3-by-2 valid tessellations that are emerged from a 2-by-2 tessellation with four TMPs shown in \textbf{a}. \textbf{c}, 2-by-3 valid tessellations that are emerged from a 2-by-2 tessellation with four TMPs shown in \textbf{a}. In both \textbf{b} and \textbf{c}, these tessellations are the only ones that can be emerged from the 2-by-2 tessellation in \textbf{a}. Cyan, red, blue, and yellow colors in the 3D rendered images correspond to the TMP, OT$+$, OT$-$, and defect states of the unit cell, respectively.} 
\label{fig:tessellationevolve}
\end{figure}

As we increase the size of the tessellation, we see that the number of unique configurations increases exponentially for $X$-by-2 and 2-by-$Y$ cases where $X$ and $Y$ are the number of horizontal and vertical stacking, respectively. For $N$-by-$N$ tessellation where $N$ is for both horizontal and vertical numbers, we see that the rate of increase decays gradually. Furthermore, we observe that the rate of increase in the number differs by the direction of adding a new column or row (horizontally or vertically). Adding a column of unit cells to the tessellation doubles the unique configurations to the original one, whereas adding a row quadruples them. In Fig.~\ref{fig:tessellationevolve}, we show an example of adding a row or a column to a 2-by-2 tessellation with four TMPs. While we obtain two different configurations in Fig.~\ref{fig:tessellationevolve}(b) for 3-by-2 tessellation, we can achieve four configurations for 2-by-3 tessellations in Fig.~\ref{fig:tessellationevolve}(c). We observe this difference is caused by the direction of the phase transition of unit cells. Since the transition among TMP, OT$+$, and OT$-$ is lateral movement, the column-wise attachment of unit cells is more bounded to the existing tessellations than the row-wise attachment. Actually, a column-wise attachment involves only TMP and OT$+$ in Fig.~\ref{fig:tessellationevolve}(b). However, a row-wise addition includes TMP, OT$+$, OT$-$, and defect in Fig.~\ref{fig:tessellationevolve}(c). 

In Fig.~\ref{fig:numberofconfig}(a), the red line connects the results from 2-by-2, 2-by-3, 2-by-4, and 2-by-8 tessellations. Likewise, the blue line connects the results from 2-by-2, 3-by-2, 4-by-2, and 8-by-2 tessellations. Both lines show linear growth in $\text{log}_2$ scale along the y-axis because there is a pattern for the growth of the unique configurations, as we state previously. The yellow line shows the results up to 7-by-7 tessellations. Since these results are from the combinations of the horizontal and vertical increments, the increase rate is not constant. 

Table.~\ref{tab:Combinations} also includes the number of unique configurations only with three states (OT$+$/OT$-$/Flat, TMP/OT$+$/Flat, TMP/OT$-$/Flat, and TMP/OT$+$/OT$-$)
and the rate of increase of the number of configurations from using three states to four states. Here, each combination of three states (OT$+$/OT$-$/Flat, TMP/OT$+$/Flat, TMP/OT$-$/Flat, and TMP/OT$+$/OT$-$) has the same number of valid configurations. The comparison between the four states indicates that each TMP state equally contributes to the excessive growth of reconfigurability. However, one of the advantages of inducing defect states is that we can achieve an enormous variety of reconfigurable shapes without incurring volume increases. Based on the graph-theoretic evaluations, we confirm that the increase of the phase number from three to four contributes directly to the surge of the overall variety of the tessellations. Further results and discussions about the effect of defect modes are included in Supplementary Note B.

\subsection{Highly heterogeneous configurations}
\begin{figure}[H]
\centering
\includegraphics[width=1.\linewidth]{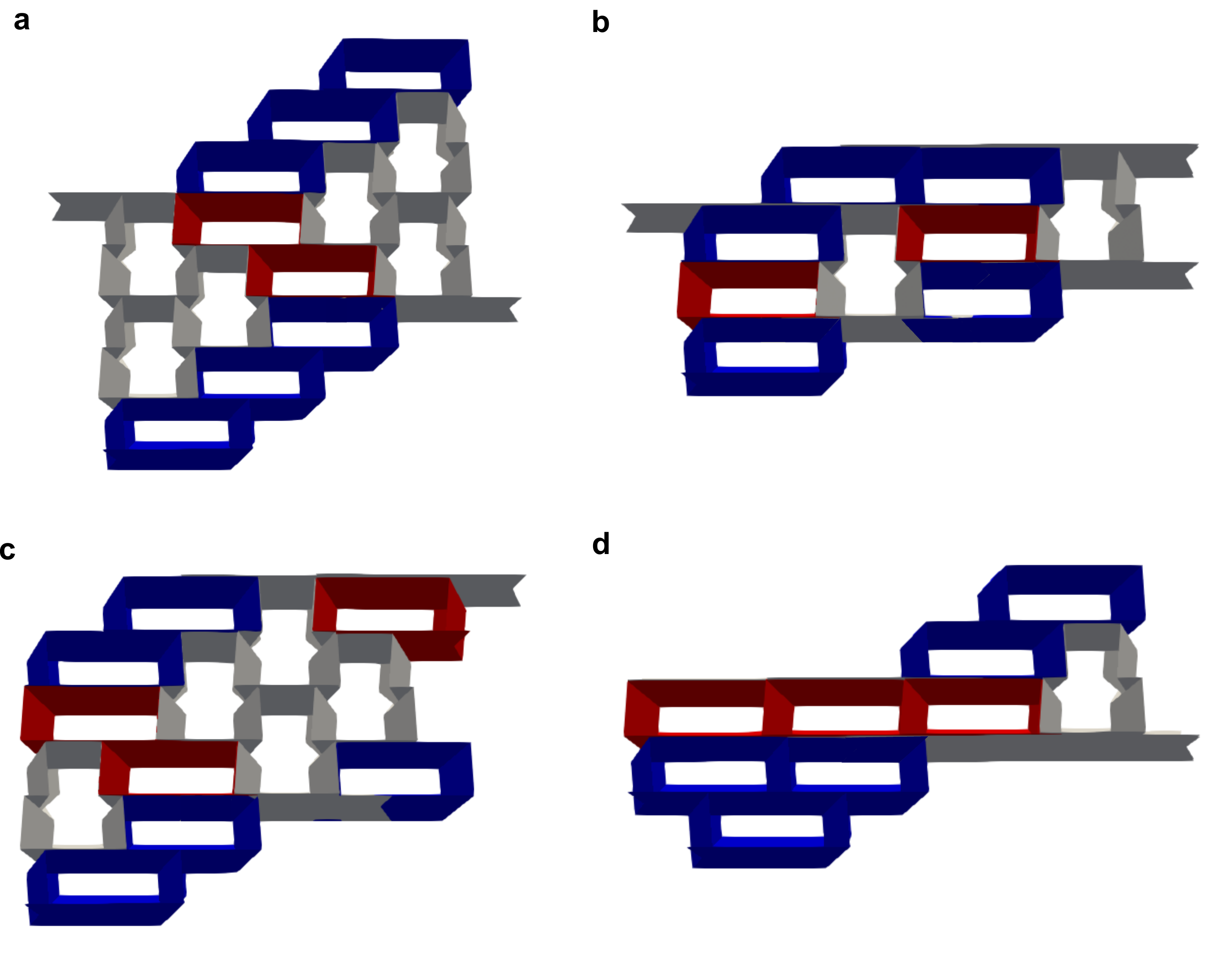}
\caption{Examples of heterogeneous configurations in 4-by-4 tessellations. \textbf{a}, a configuration consisting of 6 TMPs, 2 OT$+$, 6 OT$-$, and 2 defects. \textbf{b}, a configuration consisting of 2 TMPs, 2 OT$+$, 5 OT$-$, and 7 defects. \textbf{c}, a configuration consisting of 5 TMPs, 3 OT$+$, 5 OT$-$, and 3 defects. \textbf{d}, a configuration consisting of 1 TMP, 3 OT$+$, 5 OT$-$, and 7 defects.} 
\label{fig:anomalous}
\end{figure}

Now we discuss some examples of configurations we discover in the search process. Fig.~\ref{fig:anomalous} shows the four examples of anomalous configurations in 4-by-4 tessellations. They consist of four different states (TMP, OT$+$, OT$-$, and Defect). Also, the composition of the states is highly mixed and distributed within the tessellation. For example, a tessellation in Fig.~\ref{fig:anomalous}(a) consists of six TMPs, two OT$+$, six OT$-$, and two Defects. Besides, 3-cell groups of TMPs are separated by OT$+$, OT$-$, and defect cells. Unlike normal configurations that have a homogeneous occurrence of deformation states, such heterogeneous configurations composed of a variety of the different states are hard to be discovered by human heuristics. 

We further extend this result to see if they can simulate specific geometries because such reconfigurability is one of the essential aspects when we utilize these mechanical metamaterials for engineering applications. To assess how much the tessellation's reconfigurability contributes to the ability to simulate the real-life objects, we examine cross-sections of 4-by-4 tessellations and compare them with some images in Supplementary Note C. The result of simulating the various type of geometries (chair, buffalo, boot, and boat) suggests that this platform of TMP tessellations has a great potential to satisfy the geometrical requirements for mechanical, aerospace, and medical applications potentially.

\section{Conclusion}
In this study, we have demonstrated a graph-theoretical approach to discover the valid configurations of TMP tessellations. Given the challenges of combinatorial problems of origami-based mechanical metamaterials, we have built a computationally efficient framework to account for a complete set of reconfigurable shapes of the TMP tessellations. To this end, we have enumerated the information about the connections of the TMP unit cells using adjacency matrices. One of the unique features of this versatile graph-theoretical approach is that it has the possibility to be applied for designing various types of mechanical metamaterials. This simple approach of using adjacency matrices has proposed an efficient method to describe and analyze the tessellation of the mechanical metamaterials. The result of the aforementioned analysis indicates that the TMP tessellations have an abundance of reconfigurability owing to the heterogeneity of the TMP cells and their versatile connectivity. Such reconfigurability can be exponentially improved by adopting a larger number of TMP cells in the tessellation, and the exact number of possible configurations has been calculated by the proposed graph-theoretic method in an accurate and efficient manner.

Given the simplicity of this framework, the graph-based approach can be applied to the other types of the tessellations of mechanical metamaterials (e.g., metal-organic hinged cube tessellation~\cite{Jin2019}, voxelated mechanical metamaterials~\cite{Yang2020,Coulais2016}, and other origami lattices in 2D or 3D settings~\cite{He2020,Fang2018,Silverberg2014,Evans2015}) by building graph representations for each architecture of mechanical metamaterials and by understanding the connections within the tessellations. Also, the ability to simulate the various shapes and morphologies suggests the possibility of answering the following question; given the target geometry and properties we want to achieve, can we dial-in the local phases in unit cells and find the optimal metamaterial configurations globally? This can be assessed by utilizing the graph-theoretical framework developed herein and combining with combinatorial optimization techniques.

\section*{Declaration of competing interest}
The authors declare that they have no known competing financial interests or personal relationships that could have appeared to influence the work reported in this paper. 

\section*{CRediT authorship contribution statement}
\textbf{Koshiro Yamaguchi:} Conceptualization, Methodology, Investigation, Software, Visualization, Writing- Original draft preparation. \textbf{Hiromi Yasuda:} Conceptualization, Methodology, Investigation, Software, Visualization, Writing - Review \& Editing. \textbf{Kosei Tsujikawa:} Investigation, Software. \textbf{Takahiro Kunimine:} Supervision, Writing - Review \& Editing.  \textbf{Jinkyu Yang:} Conceptualization,  Supervision, Writing - Review \& Editing, Funding acquisition.

\section*{Acknowledgments}
K.Y. and J.Y. are grateful for the support from the U.S. National Science Foundation
(1553202 and 1933729) and the Washington Research Foundation. K.Y. is supported by the Funai Foundation for Information Technology. We thank Dr. Dillon Foight and Prof. Mehran Mesbahi at the University of Washington, and Dr. Jesse Silverberg at Multiscale Systems for helpful discussions.

\bibliographystyle{elsarticle-num}
\setlength{\bibsep}{0pt}

{\small \bibliography{origami_optim.bib}}

\begin{thebibliography}{10}
\expandafter\ifx\csname url\endcsname\relax
  \def\url#1{\texttt{#1}}\fi
\expandafter\ifx\csname urlprefix\endcsname\relax\def\urlprefix{URL }\fi
\expandafter\ifx\csname href\endcsname\relax
  \def\href#1#2{#2} \def\path#1{#1}\fi

\bibitem{Babaee2015}
S.~Babaee, N.~Viard, P.~Wang, N.~X. Fang, K.~Bertoldi, {Harnessing Deformation
  to Switch on and off the Propagation of Sound}, Advanced Materials 28~(8)
  (2016) 1631--1635.
\newblock \href {https://doi.org/10.1002/adma.201504469}
  {\path{doi:10.1002/adma.201504469}}.

\bibitem{Bilal2017}
O.~R. Bilal, A.~Foehr, C.~Daraio, {Reprogrammable Phononic Metasurfaces},
  Advanced Materials 29~(39) (2017).
\newblock \href {https://doi.org/10.1002/adma.201700628}
  {\path{doi:10.1002/adma.201700628}}.

\bibitem{Wu2018}
Z.~Wu, Y.~Zheng, K.~W. Wang, {Metastable modular metastructures for on-demand
  reconfiguration of band structures and nonreciprocal wave propagation},
  Physical Review E 97~(2) (aug 2018).
\newblock \href {http://arxiv.org/abs/1709.01800} {\path{arXiv:1709.01800}},
  \href {https://doi.org/10.1103/PhysRevE.97.022209}
  {\path{doi:10.1103/PhysRevE.97.022209}}.

\bibitem{Schenk2013}
M.~Schenk, S.~D. Guest, {Geometry of Miura-folded metamaterials}, Proceedings
  of the National Academy of Sciences of the United States of America 110~(9)
  (2013) 3276--3281.
\newblock \href {https://doi.org/10.1073/pnas.1217998110}
  {\path{doi:10.1073/pnas.1217998110}}.

\bibitem{Felton2014}
S.~Felton, M.~Tolley, E.~Demaine, D.~Rus, R.~Wood, {A method for building
  self-folding machines}, Science 345~(6197) (2014) 644--646.
\newblock \href {https://doi.org/10.1126/science.1252610}
  {\path{doi:10.1126/science.1252610}}.

\bibitem{tachi_introduction_2019}
T.~Tachi, Introduction to {Structural} {Origami}, Journal of the International
  Association for Shell and Spatial Structures 60~(1) (2019) 7--18.
\newblock \href {https://doi.org/10.20898/j.iass.2019.199.004}
  {\path{doi:10.20898/j.iass.2019.199.004}}.

\bibitem{Rus2018}
D.~Rus, M.~T. Tolley, {Design, fabrication and control of origami robots},
  Nature Reviews Materials 3~(6) (2018) 101--112.
\newblock \href {https://doi.org/10.1038/s41578-018-0009-8}
  {\path{doi:10.1038/s41578-018-0009-8}}.

\bibitem{Yang2020}
N.~Yang, C.~W. Chen, J.~Yang, J.~L. Silverberg, {Emergent reconfigurable
  mechanical metamaterial tessellations with an exponentially large number of
  discrete configurations}, Materials and Design 196 (nov 2020).
\newblock \href {https://doi.org/10.1016/j.matdes.2020.109143}
  {\path{doi:10.1016/j.matdes.2020.109143}}.

\bibitem{overvelde_three-dimensional_2016}
J.~T. Overvelde, T.~A. de~Jong, Y.~Shevchenko, S.~A. Becerra, G.~M. Whitesides,
  J.~C. Weaver, C.~Hoberman, K.~Bertoldi, A three-dimensional actuated
  origami-inspired transformable metamaterial with multiple degrees of freedom,
  Nature Communications 7~(1) (2016) 10929.
\newblock \href {https://doi.org/10.1038/ncomms10929}
  {\path{doi:10.1038/ncomms10929}}.

\bibitem{overvelde_rational_2017}
J.~T.~B. Overvelde, J.~C. Weaver, C.~Hoberman, K.~Bertoldi, Rational design of
  reconfigurable prismatic architected materials, Nature 541~(7637) (2017)
  347--352.
\newblock \href {https://doi.org/10.1038/nature20824}
  {\path{doi:10.1038/nature20824}}.

\bibitem{Filipov2015}
E.~T. Filipov, T.~Tachi, G.~H. Paulino, D.~A. Weitz, {Origami tubes assembled
  into stiff, yet reconfigurable structures and metamaterials}, Proceedings of
  the National Academy of Sciences of the United States of America 112~(40)
  (2015) 12321--12326.
\newblock \href {https://doi.org/10.1073/pnas.1509465112}
  {\path{doi:10.1073/pnas.1509465112}}.

\bibitem{yasuda_origamibased_2019}
H.~Yasuda, B.~Gopalarethinam, T.~Kunimine, T.~Tachi, J.~Yang, Origami‐{Based}
  {Cellular} {Structures} with {In} {Situ} {Transition} between {Collapsible}
  and {Load}‐{Bearing} {Configurations}, Advanced Engineering Materials
  21~(12) (2019) 1900562.
\newblock \href {https://doi.org/10.1002/adem.201900562}
  {\path{doi:10.1002/adem.201900562}}.

\bibitem{Fang2018}
H.~Fang, S.-C.~A. Chu, Y.~Xia, K.-W. Wang, {Programmable Self-Locking Origami
  Mechanical Metamaterials}, Advanced Materials 30~(15) (2018) 1706311.
\newblock \href {https://doi.org/10.1002/adma.201706311}
  {\path{doi:10.1002/adma.201706311}}.

\bibitem{Tachi2012}
T.~Tachi, K.~Miura, {Rigid-foldable cylinders and cells}, Journal of the
  International Association for Shell and Spatial Structures 53~(174) (2012)
  217--226.

\bibitem{yasuda_reentrant_2015}
H.~Yasuda, J.~Yang, Reentrant {Origami}-{Based} {Metamaterials} with {Negative}
  {Poisson}’s {Ratio} and {Bistability}, Physical Review Letters 114~(18)
  (2015) 185502.
\newblock \href {https://doi.org/10.1103/PhysRevLett.114.185502}
  {\path{doi:10.1103/PhysRevLett.114.185502}}.

\bibitem{Mesbahi2010}
M.~Mesbahi, M.~Egerstedt, Graph Theoretic Methods in Multiagent Networks.,
  Vol.~33 of Princeton Series in Applied Mathematics, Princeton University
  Press / DeGruyter, 2010.

\bibitem{Deo1974}
N.~Deo, Graph Theory with Applications to Engineering and Computer Science
  (Prentice Hall Series in Automatic Computation), Prentice-Hall, Inc., USA,
  1974.

\bibitem{Canutescu2003}
A.~A. Canutescu, A.~A. Shelenkov, R.~L. Dunbrack, {A graph-theory algorithm for
  rapid protein side-chain prediction}, Protein Science 12~(9) (2003)
  2001--2014.
\newblock \href {https://doi.org/10.1110/ps.03154503}
  {\path{doi:10.1110/ps.03154503}}.

\bibitem{Xu2018}
K.~Xu, W.~Hu, J.~Leskovec, S.~Jegelka, How powerful are graph neural networks?,
  in: International Conference on Learning Representations, 2019.

\bibitem{Arkus2011}
N.~Arkus, V.~N. Manoharan, M.~P. Brenner, {Deriving finite sphere packings},
  SIAM Journal on Discrete Mathematics 25~(4) (2011) 1860--1901.
\newblock \href {http://arxiv.org/abs/1011.5412} {\path{arXiv:1011.5412}},
  \href {https://doi.org/10.1137/100784424} {\path{doi:10.1137/100784424}}.

\bibitem{Vlassis2020}
N.~N. Vlassis, R.~Ma, W.~C. Sun, {Geometric deep learning for computational
  mechanics Part I: Anisotropic Hyperelasticity}, Computer Methods in Applied
  Mechanics and Engineering 371 (2020) 113299.
\newblock \href {http://arxiv.org/abs/2001.04292} {\path{arXiv:2001.04292}}.

\bibitem{Jin2019}
E.~Jin, I.~S. Lee, D.~Kim, H.~Lee, W.~D. Jang, M.~S. Lah, S.~K. Min, W.~Choe,
  {Metal-organic framework based on hinged cube tessellation as transformable
  mechanical metamaterial}, Science Advances 5~(5) (2019).
\newblock \href {https://doi.org/10.1126/sciadv.aav4119}
  {\path{doi:10.1126/sciadv.aav4119}}.

\bibitem{Coulais2016}
C.~Coulais, E.~Teomy, K.~{De Reus}, Y.~Shokef, M.~{Van Hecke}, {Combinatorial
  design of textured mechanical metamaterials}, Nature 535~(7613) (2016)
  529--532.
\newblock \href {http://arxiv.org/abs/1608.00625} {\path{arXiv:1608.00625}},
  \href {https://doi.org/10.1038/nature18960} {\path{doi:10.1038/nature18960}}.

\bibitem{He2020}
Y.~L. He, P.~W. Zhang, Z.~You, Z.~Q. Li, Z.~H. Wang, X.~F. Shu, {Programming
  mechanical metamaterials using origami tessellations}, Composites Science and
  Technology 189 (mar 2020).
\newblock \href {https://doi.org/10.1016/j.compscitech.2020.108015}
  {\path{doi:10.1016/j.compscitech.2020.108015}}.

\bibitem{Silverberg2014}
J.~L. Silverberg, A.~A. Evans, L.~McLeod, R.~C. Hayward, T.~Hull, C.~D.
  Santangelo, I.~Cohen, {Using origami design principles to fold reprogrammable
  mechanical metamaterials}, Science 345~(6197) (2014) 647--650.
\newblock \href {https://doi.org/10.1126/science.1252876}
  {\path{doi:10.1126/science.1252876}}.

\bibitem{Evans2015}
A.~A. Evans, J.~L. Silverberg, C.~D. Santangelo, {Lattice mechanics of origami
  tessellations}, Physical Review E 92 (2015) 13205.
\newblock \href {https://doi.org/10.1103/PhysRevE.92.013205}
  {\path{doi:10.1103/PhysRevE.92.013205}}.

\end{thebibliography}

\newpage
\section*{List of figures}

\begin{figure}[H]
\centering
\includegraphics[width=1\linewidth]{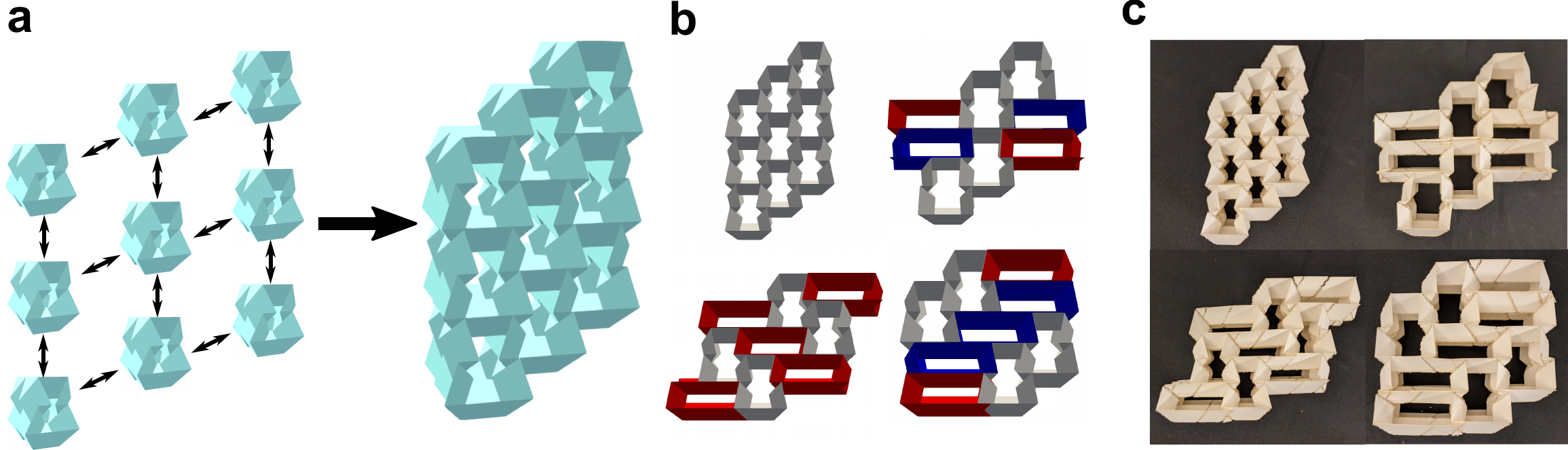}
\caption{Examples of various configurations within a 3-by-3 tessellation. \textbf{a}, 3D rendered images showing how a tessellation is built with the nine unit cells. Certain faces of unit cells are adjoined altogether to form one tessellation. \textbf{b}, 3D rendered images of the 3-by-3 tessellation with four different configurations. \textbf{c}, digital images of the paper prototypes corresponding to four tessellations shown in \textbf{b}.} 
\label{fig:tessellation}
\end{figure}

\newpage
\begin{table}[H]
    \caption{Number of unique configurations for various tessellation sizes.}
    \center
    \begin{tabular}{|l|l|l|l|l|l|l|l|}
        \hline
        Size & 4 states & 3 states & Increase & Size & 4 states & 3 states & Increase\\ \hline
        2-by-2 & 32 & 13 & 146\% & 3-by-4 & 1024 & 120 & 753\% \\ \hline
        3-by-2 & 64 & 21 & 205\% & 4-by-4 & 2048 & 184 & 1013\% \\ \hline
        2-by-3 & 128 & 32 & 300\% & 8-by-2 & 2048 & 517 & 296\% \\ \hline
        4-by-2 & 128 & 37 & 246\% & 2-by-8 & 131072 & 4181 & 3035\% \\ \hline
        2-by-4 & 512 & 89 & 475\% & 5-by-5 & 16383 & 672 & 2338\% \\ \hline
        3-by-3 & 256 & 50 & 412\% & 6-by-6 & 131072 & 2480 & 5185\% \\\hline
        4-by-3 & 512 & 82 & 524\% & 7-by-7 & 1048576 & 9312 & 11160\% \\ \hline
    \end{tabular}
    \label{tab:Combinations}
\end{table}

\end{document}


\title{\textsf{Supplementary Materials for: Graph-theoretric estimation of reconfigurability in origami-based metamaterials}}


\newpage
\maketitle

\section{\textbf{Supplementary notes}}

\subsection{TMP configurations search via brute-force method}
\begin{figure}[H]
\centering
\includegraphics[width=1.\linewidth]{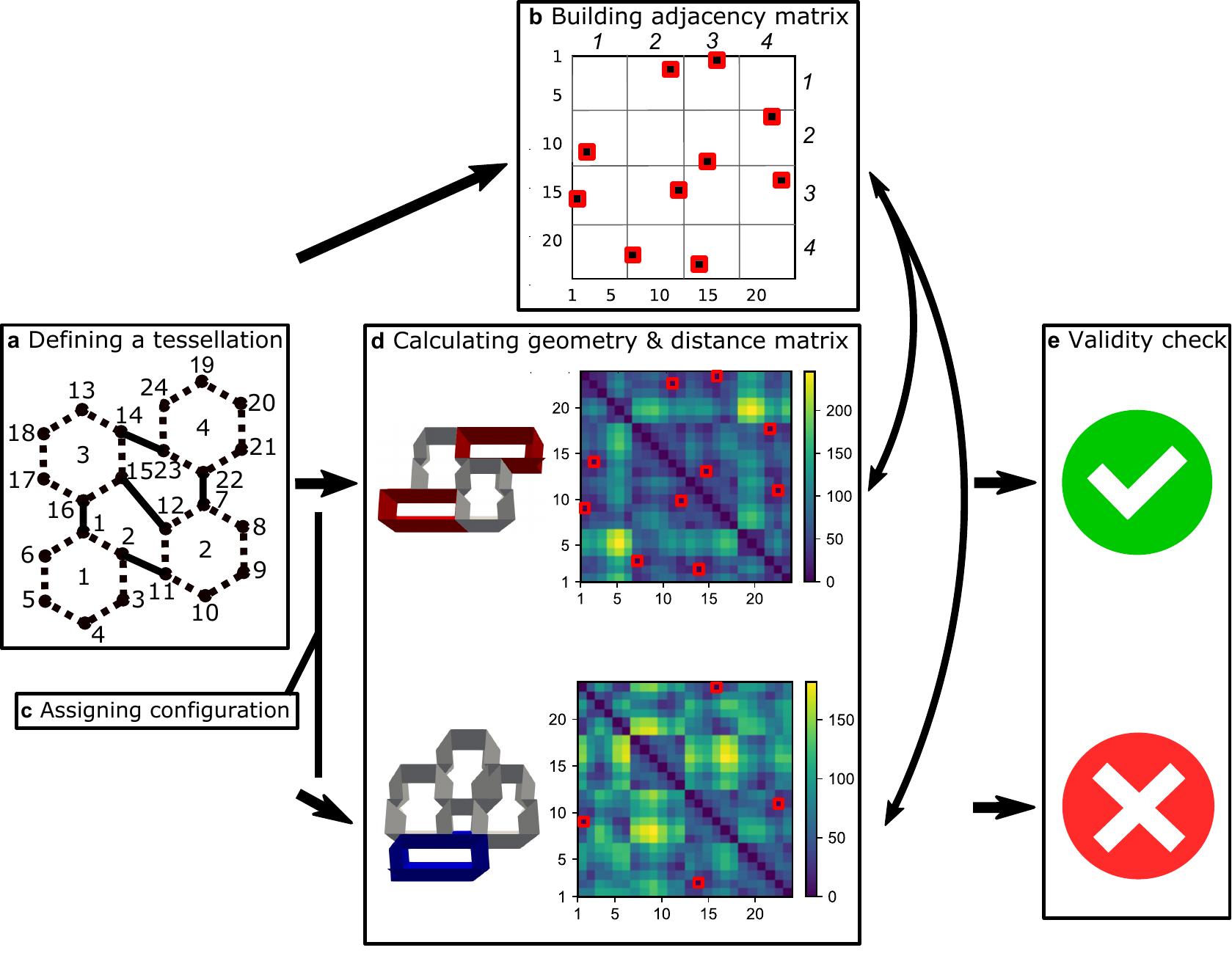}
\caption{Process of the search for valid configurations. \textbf{a}, the graph representation of 2-by-2 tessellations. \textbf{b}, a plot of a adjacency matrix. Black dots with red peripheries shows the elements of the value of 1 (nodes are connected in the subgraph $G'$). \textbf{c}, assigning configurations (TMP, RP, LP, or Flat) for each unit cell. \textbf{d}, the computation of geometry for each unit cell and the plot of distance matrices. Left images show the 3D rendering of the 2-by-2 tessellation after assigning configurations. Right images are the plot of distance matrices. Black dots with red peripheries shows the elements of the value of 0 (middle points are connected in the physical space). \textbf{e}, the result of the validity check for the configuration we assigned. If the distance matrix has entries of zero in the same position as the adjacency matrix has entries of one, we can classify a configuration as valid.}
\label{fig:searchprocess}
\end{figure}

\begin{figure}[H]
\centering
\includegraphics[width=1.\linewidth]{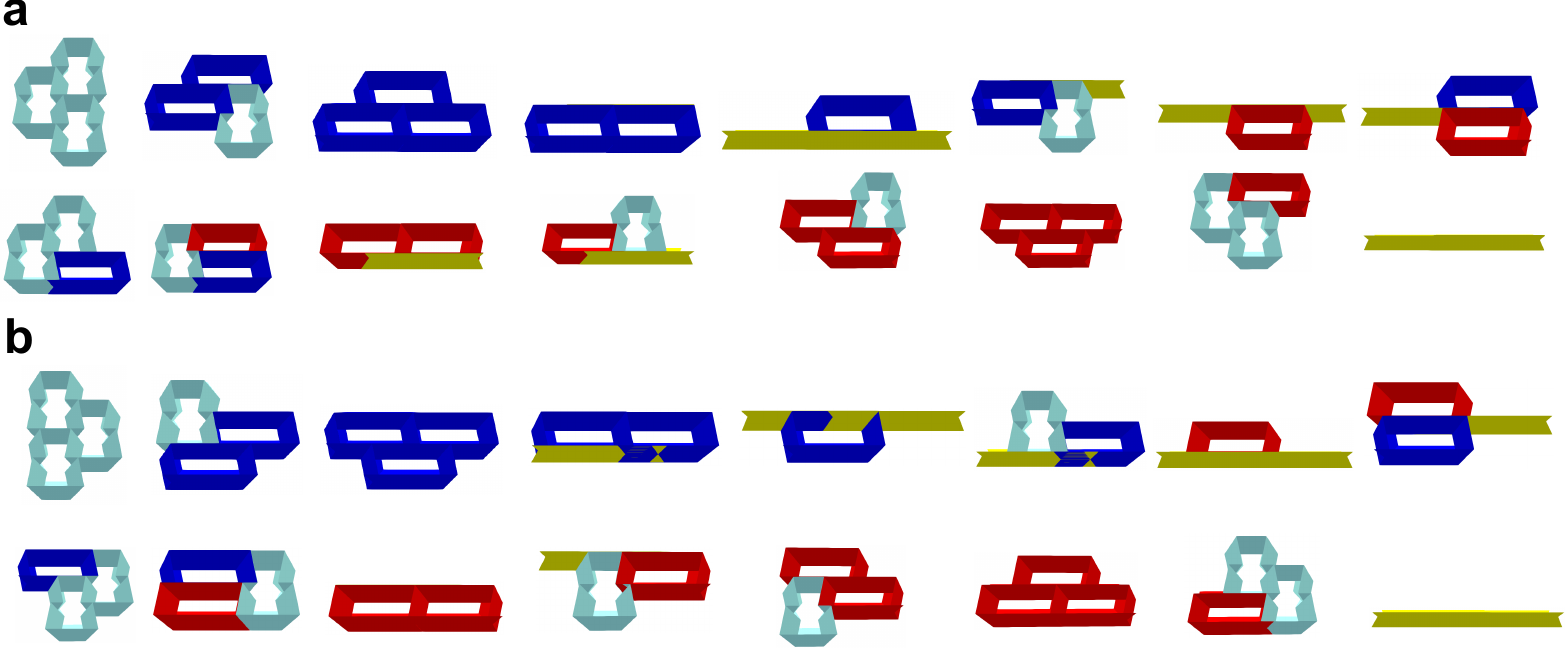}
\caption{A collection of 3-cell tessellations. \textbf{a}, tessellations that have two cells on the right side. \textbf{b}, tessellations that have two cells on the left side.}
\label{fig:3cell_library}
\end{figure}

In order to verify the result from our graph-theoretical method, we implement a brute-force method for smaller sizes of tessellations. Here, we examine the geometrical information of all $4^n$ candidates of configurations where $n$ is the number of unit cells in a tessellation. In addition to the adjacency matrix we described in the main manuscript, we utilize the distance matrix to contain the information about geometry of the tessellation and to understand how the faces of the TMPs are attached to each other. A distance matrix is a nonnegative, square, and symmetric matrix with elements corresponding to pairwise distance between the elements in a set. While several metrics can be used to build such matrix, here we refer Euclidean distance matrices as distance matrices. Here, the distance matrix is a $|V|\times |V|$ matrix $\bm{D}=(d_{ij})$. Consider a collection of $|V|$ points in a $d$-dimensional Euclidean space, assigned to the columns of matrix $\bm{X}\in \mathbb{R}^{d\times n}$, $\bm{X}=[\bm{x}_1,\bm{x}_2,...,\bm{x}_n],\bm{x}_i \in \mathbb{R}^d$, and $\bm{x}_i$ is the position vector of edge $V_i$. As we define in the main manuscript, vertices $V$ represents the six major faces of a unit cell and $|V|=6n$ where $n$ is the number of unit cells in a tessellation. Then, the squared distance $d_{ij}$ between $\bm{x}_i$ and $\bm{x}_j$ is $d_{ij} = ||\bm{x}_i - \bm{x}_j||^2 $ where $||\cdot||$ denotes the Euclidean norm. Based on the adjacency matrix and distance matrix, we can set up criteria about the validity of the configurations. The condition is $\bm{D}_{ij} = 0$ for all $i$ and $j$ where $\bm{A}_{ij} = 1$. This means that those vertices logically attached together ($\bm{A}_{ij} = 1$) have to be physically in contact ($\bm{D}_{ij} = 0$) too. The schematic illustration of this brute-force method is also shown in Fig.~\ref{fig:searchprocess}. By utilizing this method, we verify that the number of valid configuration up to 4-by-4 tessellations matches the results from the efficient method in the main manuscript. Also, this brute-force method verifies the list of 3-cell tessellations in Fig.~\ref{fig:3cell_library}. These 3-cell tessellations are the fundamental building blocks for any size of TMP tessellations and therefore it is important to verify all of them by using this method.

\subsection{Analysis on the effect of defect modes}
\begin{figure}[H]
\centering
\includegraphics{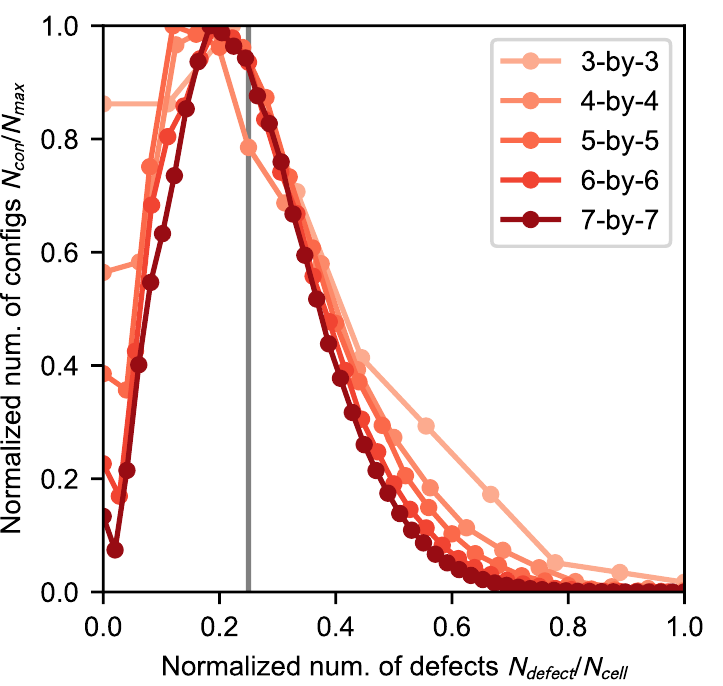}
\caption{Normalized number of configurations to normalized number of defects. Grey vertical line indicates the value of 0.25 in the x-axis. Five series of data from $n$-by-$n$ tessellation has different number of data points with connecting lines.} 
\label{fig:numberofdefects}
\end{figure}

As we mention in the main manuscript, the induction of non-voluminous defects states are not usual idea to build tessellations. However, we find that this is highly advantageous to greatly increase the reconfigurability. Here, we discuss the effect of inducing defect modes on TMP tessellations.
Fig.~\ref{fig:numberofdefects} shows the normalized number of configurations $N_{con}/N_{max}$ with the normalized number of defects $N_{defect}/N_{cell}$ for each number of defects. Here, $N_{con}$ is the number of valid configurations discovered for each number of defects. $N_{max}$ is the maximum number of configurations at a certain number of configurations. $N_{defect}$ is the number of defects in a tessellation. $N_{cell}$ is the number of unit cells in a tessellation. In Fig.~\ref{fig:numberofdefects}, the peak of the number of configurations moves rightward as the size of the tessellation enlarges. This means that we need a larger ratio of defect cells in the tessellation to get the maximum normalized number of configurations to the normalized number of defects as we have a larger size of tessellations. As the size of tessellation gets larger, we observe that the peak of the normalized number of defects makes the transition to the rightward. This implies that we need more proportion of defects to obtain maximum variation of reconfigurability. We estimate this peak eventually limits to the one-forth of the normalized number of unit cells as we have larger tessellations. The investigation on the proportion of each state to obtain maximum reconfigurability for larger tessellations is of future interest.

\begin{figure}[H]
\centering
\includegraphics{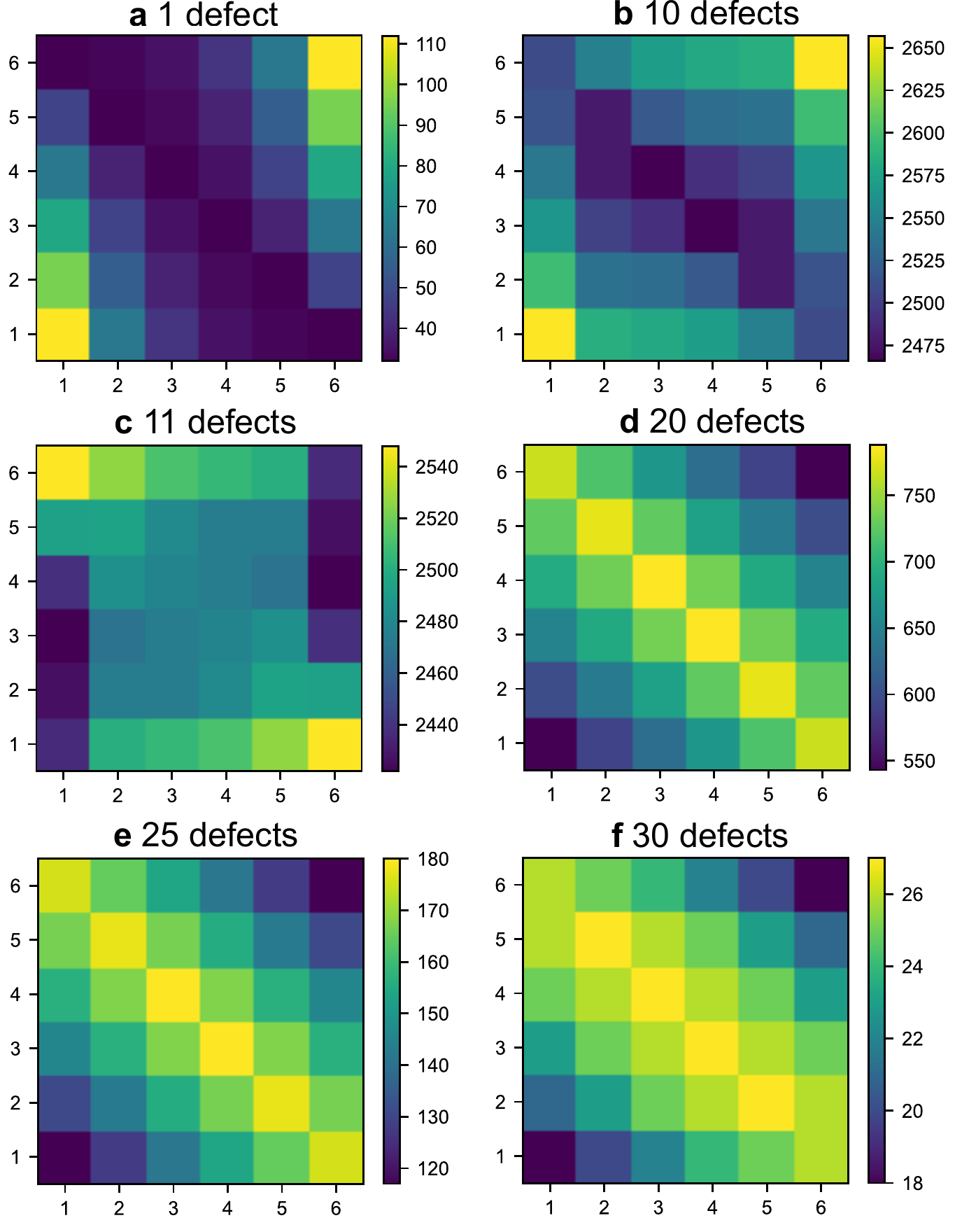}
\caption{The occurrence of defects to the given number of defects in a 6-by-6 tessellation. The numbers in $x$-axis and $y$-axis denote the position of the unit cell. Color bars for each figure shows the number of occurrence. Yellow color in the figure shows that there is a high occurrence of defects whereas blue color shows a low occurrence.} 
\label{fig:defectsposition}
\end{figure}

Furthermore, we find an intriguing transition of pattern of the occurrence of defect cells by changing the number of defect to embed in tessellations. Fig.~\ref{fig:defectsposition} shows the number of the occurrence of defect cells under a given number of defect cells in a tessellation. For instance, Fig.~\ref{fig:defectsposition}b indicates that 2650 configurations have a defect at the position of $(1,1)$ cell in case of 10-defect tessellations, whereas there are 2475 configurations that have a defect at $(3,4)$ and $(4,3)$ position. Fig.~\ref{fig:defectsposition}a to f show that there is a transition of positions where we can find the maximum number of the occurrence of the defect cell. In Fig.~\ref{fig:defectsposition}a and Fig.~\ref{fig:defectsposition}b, bottom left $(1,1)$ and top right $(6,6)$ have the largest number of defect occurrence. Then, in Fig.~\ref{fig:defectsposition}c, top left $(1,6)$ and bottom right $(6,1)$ have the largest occurrence. Afterwards, the peak of occurrence makes a transition to the center of the tessellations as shown in Fig.~\ref{fig:defectsposition}d, e, and f. These results suggests that we can embed any number of defect cells in a tessellation. However, the position and the number of defects greatly affects the reconfigurability of the tessellation.

\subsection{Cross-sectional similarity to real-life objects}
\begin{figure}[H]
\centering
\includegraphics[width=1.\linewidth]{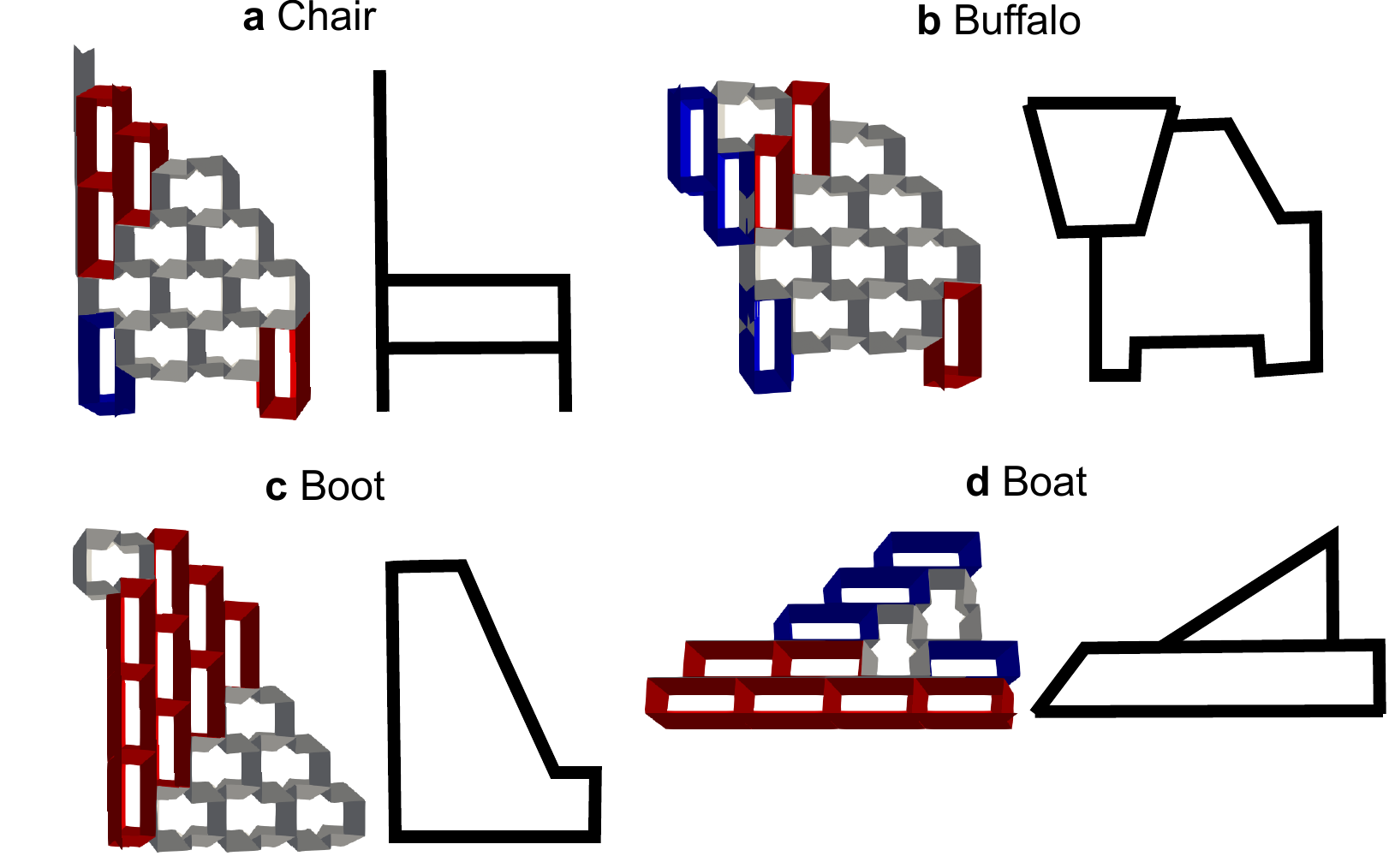}
\caption{Voxel-art configs in 4by4. In \textbf{a}, \textbf{b}, \textbf{c}, and \textbf{d}, each picture of 4-by-4 tessellations has corresponding images of chair, buffalo, boot, and boat, respectively. Pictures in \textbf{a}, \textbf{b}, and \textbf{c} are rotated by 90 degrees.} 
\label{fig:voxelart}
\end{figure}

In the main manuscript, we find that TMP tessellations have reconfigurability to a great extent. Here, we introduce some of the configurations whose cross sections have similarity to real-life objects. Fig.~\ref{fig:voxelart} shows voxel-art-like cross sections of 4-by-4 TMP tessellations that have resemblance with a chair, buffalo, boot, and boat. While they only have 16 unit cells, these results imply that the rich reconfigurability of TMP tessellation can be engineered to meet the complex geometric design requirements for the actual engineering use with increasing the number of unit cells in a tessellation.

\subsection{Measuring shape resemblance}
\begin{figure}[H]
\centering
\includegraphics[width=1.\linewidth]{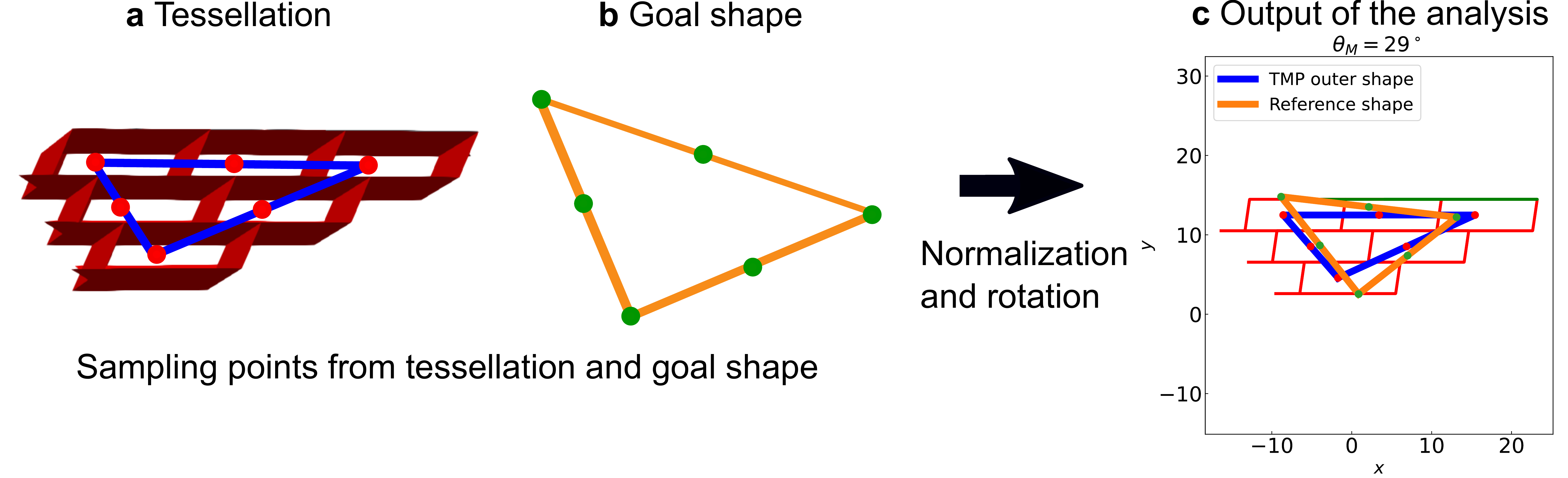}
\caption{Process of Procrustes analysis for a TMP tessellation. In \textbf{a}, red dots represent sampling points from the geometry of 3-by-3 TMP tessellation composed of five OT+ and four flat units. Blue lines connect the red dots and they represent the characteristic shape for a TMP tessellation. In \textbf{b}, orange lines represent our reference shape to be compared with the TMP tessellation. Green dots are sampling points from the reference shape. In \textbf{c}, the plot shows the result of the shape-matching calculation after applying normalization and rotation. Red and green lines shows the original geometry of the TMP. The sum of the pairwise distance between each red dot and green dot becomes the measurement of shape similarity.} 
\label{fig:procrustes}
\end{figure}

Procrustes analysis [1] is a test to check the similarity between two data sets. It is a rigid shape analysis that utilizes scaling, translation, and rotation to minimize the sum of the distance between two points in each data set. Fig.~\ref{fig:procrustes} illustrates how we sample landmarks as data sets from both the 2D cross-section of a tessellation and a goal shape. With combining the ideas in Fig.~\ref{fig:voxelart}, this method can be applied for the further investigation to search the configurations with the desired cross section that meets the engineering requirements. The details of the method is explained as follows: we sample the landmarks from the center points of the perimetral non-flat unit cells, where we define the same number of landmarks along with the goal shape. Then, we define two matrices $\bm{A} = \{\bm{a}_1, ..., \bm{a}_n\}$ and $\bm{B} = \{\bm{b}_1, ..., \bm{b}_n\}$, where $\bm{a}_i \in \mathbb{R}^2$ and $\bm{b}_i \in \mathbb{R}^2$. The generalized procrustes analysis first centers two data sets to the origin by subtracting mean values $\Bar{\bm{A}}$ and $\Bar{\bm{B}}$.

\begin{equation}
    \begin{split}
        &\bm{A}_o = \bm{A} - \Bar{\bm{A}}\\
        &\bm{B}_o = \bm{B} - \Bar{\bm{B}}
    \end{split}
\end{equation}

Then we normalize the scales of two data sets by dividing them by their own Frobenius norms.

\begin{equation}
    \begin{split}
        &\bm{A}_\text{norm} = \frac{\bm{A}_o}{||\bm{A}_o||_F}\\
        &\bm{B}_\text{norm} = \frac{\bm{B}_o}{||\bm{B}_o||_F}
    \end{split}
\end{equation}

After redefining $\bm{A}_\text{norm}$ and $\bm{B}_\text{norm}$ as $\bm{A}$ and $\bm{B}$, we can obtain the rotation matrix $\bm{R}$ that minimizes the error between $\bm{A}$ rotated by $\bm{R}$ and $\bm{B}$ via the following equation:
\begin{equation}
    \bm{R} = \text{argmin}_{\bm{\Omega}}||\bm{\Omega}\bm{A}-\bm{B}||_F \quad \text{s.t.} \quad \bm{\Omega}^T\bm{\Omega}=\bm{I}
\end{equation}

Here,
\begin{equation}
    ||\bm{\Omega}\bm{A}-\bm{B}||^2_F=||\bm{A}||_F^2+||\bm{B}||_F^2-2\langle \bm{\Omega}\bm{A},\bm{B}\rangle
\end{equation}

Therefore, the minimization of this term means maximization of the third term of the right side. Then, with singular value decomposition so that we get $\bm{A}\bm{B}^T=\bm{U}\bm{\Sigma}\bm{V}^T$:
\begin{equation}
    \begin{split}
        \langle \bm{\Omega}\bm{A}, \bm{B} \rangle &= \text{tr}(\bm{A}^T\bm{\Omega}^T\bm{B})=\text{tr}(\bm{B}^T\bm{\Omega}\bm{A})=\text{tr}(\bm{\Omega}\bm{A}\bm{B}^T)\\
        &=\text{tr}(\bm{\Omega}\bm{U}\bm{\Sigma}\bm{V}^T)=\text{tr}(\bm{\Sigma}\bm{V}^T\bm{\Omega}\bm{U})
    \end{split}
\end{equation}

If we define $\bm{X}=\bm{V}^T\bm{\Omega}\bm{U}$, $\bm{X}$ is orthogonal. Therefore,
\begin{equation}
    \text{tr}(\bm{\Sigma}\bm{V}^T\bm{\Omega}\bm{U})=\text{tr}(\bm{\Sigma}\bm{X})=\Sigma_i s_{ii}x_{ii}
\end{equation}
$s_{ii}$ is a singular value and nonnegative, therefore $\bm{X}=\bm{I}$ maximizes $\text{tr}(\bm{SX})$. Thus, $\bm{\Omega}=\bm{V}\bm{U}^T$ minimizes the original equation and that $\bm{\Omega}$ is the rotation matrix $\bm{R}$ as a solution.

\subsection*{Reference}
\noindent
[1] W. J. Krzanowski, Principles of Multivariate Analysis: A User’s Perspective, Oxford University Press, 1988 
\clearpage

\section{Supplementary figures}

